\documentclass[11pt]{article}
\usepackage{latexsym,amssymb,wrapfig}
\usepackage[all]{xy}
\def\demo{\noindent{\bf Proof. }}
\def\QED{\hfill$\Box$}
\newtheorem{Theorem}{Theorem}[section]

\newtheorem{Corollary}[Theorem]{Corollary}
\newtheorem{Proposition}[Theorem]{Proposition}
\newtheorem{Example}[Theorem]{Example}
\newtheorem{Definition}[Theorem]{Definition}

\begin{document}
\topmargin3mm
\hoffset=-1cm
\voffset=-1.5cm
\

\medskip

\begin{center}
{\large\bf Systems with the integer rounding property in normal
monomial subrings 
}
\vspace{6mm}\\
\footnotetext{2000 {\it Mathematics Subject 
Classification}. Primary 13B22; Secondary 13H10;13F20;52B20.} 
\footnotetext{{\it Key words and phrases\/}.  
canonical module, a-invariant,normal ideal, perfect graph, maximal
cliques, Rees algebra, Ehrhart ring, integer rounding property}

\medskip

Luis A. Dupont, Rafael H. Villarreal
\footnote{Partially supported by CONACyT 
grant 49251-F and SNI.} 
\\ 
{\small Departamento de Matem\'aticas}\vspace{-1mm}\\ 
{\small Centro de Investigaci\'on y de Estudios Avanzados del
IPN}\vspace{-1mm}\\   
{\small Apartado Postal 14--740}\vspace{-1mm}\\ 
{\small 07000 M\'exico City, D.F.}\vspace{-1mm}\\ 
{\small e-mail: {\tt vila@math.cinvestav.mx}}\vspace{4mm}

Carlos Renter\'\i a-M\'arquez\footnote{Partially supported by
COFAA-IPN and SNI.} 
\\
{\small Departamento de Matem\'aticas}\vspace{-1mm}\\ 
{\small Escuela Superior de F\'\i sica y
Matem\'aticas}\vspace{-1mm}\\  
{\small Instituto Polit\'ecnico Nacional}\vspace{-1mm}\\   
{\small 07300 M\'exico City, D.F.}\vspace{-1mm}\\  
{\small e-mail: {\tt renteri@esfm.ipn.mx}}\vspace{-1mm}\\  

\end{center}
\date{}

\begin{abstract}\noindent 
Let $\mathcal{C}$ be a clutter and let $A$ be its incidence matrix. If
the linear system $x\geq 0;xA\leq\mathbf{1}$ has the integer rounding
property, we give a description of the canonical module and the
$a$-invariant of certain normal subrings associated to $\mathcal{C}$.
If the clutter is a connected graph, we 
describe when the aforementioned linear system has the integer
rounding property 
in combinatorial and algebraic 
terms using graph theory and the theory of Rees algebras. As a
consequence we show that the extended Rees algebra of 
the edge ideal of a bipartite graph is Gorenstein if and only if the
graph is unmixed.
\end{abstract}

\section{Introduction}\label{intro-sept8-09}

A {\it clutter\/} $\mathcal{C}$ with finite vertex set $X$ is a 
family of subsets of $X$, 
called edges, none of which is included in another. The set of
vertices and edges of $\mathcal{C}$ are denoted by $V(\mathcal{C})$
and $E(\mathcal{C})$ respectively. A basic example  of a clutter is a
graph.  

Let $\mathcal{C}$ be a clutter with finite vertex set 
$X=\{x_1,\ldots,x_n\}$. We shall always assume that 
$\mathcal{C}$ has no isolated vertices, i.e., each vertex occurs in
at least one edge. Let $f_1,\ldots,f_q$ be the edges of $\mathcal{C}$
and 
let $v_k=\sum_{x_i\in f_k}e_i$ be the
{\it characteristic vector\/} of $f_k$, where $e_i$ is the 
$i${\it th} unit vector in $\mathbb{R}^n$. The {\it incidence matrix\/}
$A$ of $\mathcal{C}$ is the $n\times q$ matrix with column vectors
$v_1,\ldots,v_q$. If $a=(a_i)$ and $c=(c_i)$ are vectors 
in $\mathbb{R}^n$, then $a\leq c$ means that $a_i\leq c_i$ for 
all $i$. Thus $a\geq 0$, means that $a_i\geq 0$ for all $i$.
The system $x\geq 0; xA\leq\mathbf{1}$  
has the {\it integer rounding property\/} if 
$$\lceil{\rm min}\{\langle y,{\mathbf 1}\rangle \vert\, 
y\geq 0;\, Ay\geq \alpha \}\rceil
={\rm min}\{\langle y,{\mathbf 1}\rangle \vert\, 
Ay\geq \alpha;\, y\in\mathbb{N}^q\}
$$
for each integral vector $\alpha$ for which 
${\rm min}\{\langle y,{\mathbf 1}\rangle \vert\, 
y\geq 0;\, Ay\geq \alpha\}$ is finite. Here $\mathbf{1}$ denotes the
vector with all its entries equal to $1$ and $\langle\ ,\, \rangle$ 
denotes the standard 
inner product. For a thorough study of this property see (Schrijver
1986).  

Let $R=K[x_1,\ldots,x_n]$ be a polynomial ring over a field $K$ and
let $w_1,\ldots,w_r$ be the set of all integral vectors $\alpha$ such
that $0\leq \alpha\leq v_i$ for some $i$. We will examine the integer
rounding property using the monomial subring:
$$
S=K[x^{w_1}t,\ldots,x^{w_r}t]\subset R[t],
$$
where $t$ is a new variable. As usual we
use the notation $x^a:=x_1^{a_1} \cdots x_n^{a_n}$, 
where $a=(a_1,\ldots,a_n)\in \mathbb{N}^n$. The subring $S$ is called
{\it normal\/} if $S$ is integrally closed, i.e., $S=\overline{S}$, 
where $\overline{S}$ is the integral closure of $S$ in its field of
fractions, see (Vasconcelos 2005). 

The contents of this paper are as follows. One of the results 
in (Brennan et al. 2008) shows that  
$S$ is normal if and only if the system $x\geq 0;\, xA\leq
\mathbf{1}$ has 
the integer rounding property (see Theorem~\ref{round-up-char}). 
As a consequence we show that if all edges of $\mathcal{C}$ have the
same number of elements and either linear system $x\geq 0; xA\leq\mathbf{1}$
or $x\geq 0; xA\geq\mathbf{1}$ has the integer 
rounding property, then the subring $K[x^{v_1}t,\ldots,x^{v_q}t]$ is
normal (see Corollary~\ref{march6-08}). 

Let $G$ be a connected graph and let $A$ be its
incidence matrix. The main results of
Section~\ref{rounding-incidence} show that the following conditions 
are equivalent (see
Theorems~\ref{round-for-conn-graphs} and 
\ref{roundup-iff-rounddown}):  
\begin{description} 
\item{\rm (a)} $x\geq 0 ;xA\leq\mathbf{1}$ has the integer rounding
property.
\item{\rm (b)} $x\geq 0;xA\geq\mathbf{1}$ has the integer rounding
property (see Definition~\ref{irp>=}). 

\item{\rm (c)} $R[x^{v_1}t,\ldots,x^{v_q}t]$ is normal, where
$v_1,\ldots,v_q$ are the column 
vectors of $A$. 
\item{\rm (d)} $K[x^{v_1}t,\ldots,x^{v_q}t]$ is normal.

\item{\rm (e)} $K[t,x_1t,\ldots,x_nt,x^{v_1}t,\ldots,x^{v_q}t]$ is
normal.
\item{\rm (f)} The induced 
subgraph of the vertices of any two vertex disjoint odd cycles of $G$ is
connected.
\end{description}

The most interesting part of this result is the equivalence between
(a), which is a linear optimization property, and (f), which
is a graph theoretical property. Edge ideals are defined in
Section~\ref{section-introunding-normality}. We prove that the
ring in (e) 
is isomorphic to the 
extended Rees algebra of the edge ideal of $G$ 
(see Proposition~\ref{ext-rees-alg=monsub}). If $G$ is bipartite and
$I=I(G)$ 
is its edge ideal, we are able to prove that the extended Rees
algebra of $I$ is a Gorenstein standard $K$-algebra if and only if $G$
is unmixed (see Corollary~\ref{extended-gorenstein}). If we work in
the more general context of clutters none of the conditions 
(a) to (e) are equivalent. Some of these conditions are equivalent
under certain assumptions (Gitler et al. 2009). 

In Section~\ref{section-on-canmod} we introduce the canonical module
and the $a$-invariant of $S$. This invariant plays a key role in the
theory of Hilbert functions (Bruns and Herzog 1997). There are some
methods, based on combinatorial optimization, 
that have been used to study canonical modules of edge
subrings of bipartite graphs (Valencia and Villarreal 2003). Our
approach to study 
canonical modules is inspired by these methods. If
$S$ is a normal domain, we express the canonical module of $S$ and its
$a$-invariant 
in terms of the vertices of the polytope 
$$
\{x\vert x\geq
0;xA\leq\mathbf{1}\}
$$ 
(see Theorem~\ref{can-mod-intr}). We are able to give an explicit
description of 
the canonical module of $S$ and
its $a$-invariant when $\mathcal{C}$ is the clutter of maximal cliques
of a perfect graph (Theorem~\ref{perfect-canon-ainv}). 

For unexplained terminology and 
notation on commutative algebra and integer programming we refer to
(Bruns and Herzog 1997, Vasconcelos 2005, Villarreal 2001) and
(Schrijver 1986) respectively.

\section{Integer rounding and
normality}\label{section-introunding-normality}

We continue
using the 
definitions and terms from the introduction. In what follows
$\mathbb{N}$ 
denotes the set of non-negative
integers and $\mathbb{R}_+$ denotes the set of non-negative real
numbers. Let $\mathcal{A}\subset\mathbb{Z}^n$. The convex hull of
$\mathcal{A}$ is denoted by ${\rm conv}(\mathcal{A})$ and the
cone generated by $\mathcal{A}$ is denoted by
$\mathbb{R}_+\mathcal{A}$. 

\begin{Theorem}[\rm Brennan et al. 2008]\label{round-up-char} Let
$\mathcal{C}$ be a clutter 
and let 
$v_1,\ldots,v_q$ be the columns of the incidence matrix $A$ of
$\mathcal{C}$. If $w_1,\ldots,w_r$ is the set of all
$\alpha\in\mathbb{N}^n$ such that $\alpha\leq v_i$ for some $i$, 
then the system $x\geq 0;\, xA\leq \mathbf{1}$ has the integer
rounding 
property if and only if the subring
$K[x^{w_1}t,\ldots,x^{w_r}t]$ is normal. 
\end{Theorem}

Next we give an application of this result, but first we need to
introduce some more 
terminology and notation. We have already defined in the introduction
when the linear  
system $x\geq 0; xA\leq\mathbf{1}$ has the integer rounding property.
The following is a dual notion.

\begin{Definition}\label{irp>=}\rm Let $A$ be a matrix with entries in
$\mathbb{N}$. The system $x\geq 0; xA\geq\mathbf{1}$ 
has the {\it integer rounding property\/} if 
$${\rm max}\{\langle y,{\mathbf 1}\rangle \vert\, 
Ay\leq w; y\in\mathbb{N}^q\} 
=\lfloor{\rm max}\{\langle y,{\mathbf 1}\rangle \vert\, y\geq 0;
Ay\leq w\}\rfloor
$$
for each integral vector $w$ for which the right 
hand side is finite. 
\end{Definition}

Let $\mathcal{A}=\{v_1,\ldots,v_q\}$ be a
set 
of points in $\mathbb{N}^n$, let 
$P$ be the convex hull of $\mathcal{A}$, and 
let $R=K[x_1,\ldots,x_n]$ be a polynomial ring 
over a field $K$. The {\it Ehrhart ring\/} of
the lattice polytope $P$ is the monomial subring 
$$
A(P)=K[\{x^at^i\vert\, a\in \mathbb{Z}^n \cap iP; i\in
\mathbb{N}\}]\subset R[t],
$$
where $t$ is a new variable. 
A nice property of $A(P)$ is that it is always a normal domain 
(Bruns and Herzog 1997). Let $\mathcal{C}$ be a clutter with vertex
set 
$X=\{x_1,\ldots,x_n\}$ 
and edge set $E(\mathcal{C})$. 
For use below recall that $\mathcal{C}$ is called {\it
uniform\/} if all its edges have the same number of elements. The
{\it edge 
ideal\/} of $\mathcal{C}$,  
denoted by $I(\mathcal{C})$, is the ideal of $R$
generated by all monomials $\prod_{x_i\in e}x_i=x_e$ such 
that $e\in E(\mathcal{C})$. The {\it Rees algebra\/} of $I=I(\mathcal{C})$, 
denoted by $R[It]$, is given by 
$$
R[It]:=R\oplus It\oplus\cdots\oplus I^{i}t^i\oplus\cdots
\subset R[t],
$$
see (Vasconcelos 2005, Chapter~1) for a nice presentation of Rees
algebras. 

\medskip

The next result holds for arbitrary monomial ideals (not necessarily
square-free). This is the only 
place in the paper where a result is stated for arbitrary monomial
ideals.  

\begin{Theorem}[\rm Dupont and Villarreal
2010]\label{poset-main-normal} Let 
$I=(x^{v_1},\ldots,x^{v_q})\subset R$ be a monomial ideal
and let $A$ be the matrix with column vectors $v_1,\ldots,v_q$. Then
the system $x\geq 0; xA\geq\mathbf{1}$ has the integer rounding
property if and only if $R[It]$ is normal.
\end{Theorem}

\begin{Corollary}\label{march6-08} Let $\mathcal{C}$ be a uniform
clutter, let $A$ 
be its incidence matrix, and let $v_1,\ldots,v_q$ be the columns of
$A$. If either system $x\geq 0; xA\leq\mathbf{1}$
or $x\geq 0; xA\geq\mathbf{1}$ has the integer 
rounding property and $P={\rm conv}(v_1,\ldots,v_q)$, then 
$$K[x^{v_1}t,\ldots,x^{v_q}t]=A(P).$$
\end{Corollary}

\demo  In general the subring $K[x^{v_1}t,\ldots,x^{v_q}t]$ is
contained in $A(P)$.  
Assume that $x\geq 0; xA\leq\mathbf{1}$ has the integer rounding
property and that every edge of $\mathcal{C}$ has $d$ 
elements. 
Let $w_1,\ldots,w_r$ be the set of all
$\alpha\in\mathbb{N}^n$ such that $\alpha\leq v_i$ for some $i$. Then
by Theorem~ \ref{round-up-char} the subring
$K[x^{w_1}t,\ldots,x^{w_r}t]$ is normal. Using that $v_1,\ldots,v_q$
is the set of $w_i$ with $|w_i|=d$, it is not hard to see that 
$A(P)$ is contained in $K[x^{v_1}t,\ldots,x^{v_q}t]$. 

Assume that $x\geq 0; xA\geq\mathbf{1}$ has the integer rounding
property. Let $I=I(\mathcal{C})$ be the edge ideal of $\mathcal{C}$
and let $R[It]$ be its Rees algebra. By
Theorem~\ref{poset-main-normal}, 
$R[It]$ is a
normal domain. Since
the clutter $\mathcal{C}$ is uniform the required equality 
follows at once from (Escobar et al. 2003, Theorem~3.15). \QED

\medskip

The converse of Corollary~\ref{march6-08} fails as the following
example shows. 

\begin{Example}\rm Let $\mathcal{C}$ be the uniform clutter with 
vertex set $X=\{x_1,\ldots,x_8\}$ and edge set
$$
E(\mathcal{C})=\{\{x_3,x_4,x_6,x_8\},\, \{x_2,x_5,x_6,x_7\},\,
\{x_1,x_4,x_5,x_8\},\, \{x_1,x_2,x_3,x_8\}\}. 
$$
The characteristic vectors of the edges of $\mathcal{C}$ are
\begin{eqnarray*}
v_1=(0, 0, 1, 1, 0, 1, 0, 1),&v_2=(0, 1, 0, 0, 1, 1, 1, 0),&\\
v_3=(1, 0, 0, 1, 1, 0, 0, 1),&v_4=(1, 1, 1, 0, 0, 0, 0, 1).&
\end{eqnarray*}
Let $A$ be the incidence matrix of $\mathcal{C}$ with column vectors
$v_1,\ldots,v_4$ and let $P$ be the convex hull of
$\{v_1,\ldots,v_4\}$. It is not 
hard to verify that the set 
$$\{(v_1,1),(v_2,1),(v_3,1),(v_4,1)\}$$ 
is a Hilbert basis in the sense of (Schrijver 1986). Therefore we
have the 
equality
$$K[x^{v_1}t,x^{v_2}t,x^{v_3}t,x^{v_4}t]=A(P).$$
Using Theorem~\ref{poset-main-normal} and (Brennan et al. 2008,
Theorem~2.12) it is seen that none of the two systems $x\geq 0;
xA\leq\mathbf{1}$ 
and $x\geq 0; xA\geq\mathbf{1}$ have the integer 
rounding property.  
\end{Example}

\section{Incidence matrices of graphs}
\label{rounding-incidence}

Let $G$ be a connected graph with vertex set $X=\{x_1,\ldots,x_n\}$
and let $v_1,\ldots,v_q$ be the column vectors of the incidence 
matrix $A$ of $G$. The main result here is a purely combinatorial
description 
of the integer rounding property of the system $x\geq 0;xA\leq
\mathbf{1}$.
Other equivalent algebraic conditions of this property will be
presented.  

Let $R=K[x_1,\ldots,x_n]$ be a polynomial ring over
a field $K$ and let $I=I(G)$ be the edge ideal of $G$. Recall that the

{\it extended Rees algebra\/} of $I$ is the subring 
$$
R[It,t^{-1}]:=R[It][t^{-1}]\subset R[t,t^{-1}],
$$
where $R[It]$ is the Rees algebra of $I$. Rees algebras of edge 
ideals of graphs were first studied in (Simis et al. 1994).

\begin{Proposition}\label{ext-rees-alg=monsub}
$R[It,t^{-1}]\simeq K[t,x_1t,\ldots,x_nt,x^{v_1}t,\ldots,x^{v_q}t]$. 
\end{Proposition}

\demo  We set $S=K[t,x_1t,\ldots,x_nt,x^{v_1}t,\ldots,x^{v_q}t]$. Note
that $S$ and $R[It,t^{-1}]$ are both integral domains of the same
Krull dimension, this follows from the dimension formula given in 
(Sturmfels 1996, Lemma~4.2). Thus it suffices to prove that there is an 
epimorphism $\overline{\psi}\colon S\rightarrow R[It,I^{-1}]$ of
$K$-algebras.

Let $u_0,u_1,\ldots,u_n,t_1,\ldots,t_q$ be a
new set of variables and let $\varphi$, $\psi$ be the maps of
$K$-algebras defined  
by the diagram
\vspace{-8mm}
$$
\begin{array}{ccc}
\xymatrix{
K[u_0,u_1,\ldots,u_n,t_1,\ldots,t_q]\ar[r]^{\ \ \ \ 
 \ \ \ \ \ \ \psi} \ar[d]^\varphi&R[It,t^{-1}]
& \ \ \\ 
S
\ar@{-->}[ru]^{\overline{\psi}}&
}\hspace{-7mm}
&
\begin{array}{l}
\ \\
\ \\
u_0\stackrel{\varphi}{\longmapsto} t,\\ 
u_i\stackrel{\varphi}{\longmapsto} x_it,\\
t_i\stackrel{\varphi}{\longmapsto} x^{v_i}t,
\end{array}
&
\begin{array}{l}
\ \\ 
\ \\
u_0\stackrel{\psi}{\longmapsto} t^{-1}, \\
u_i\stackrel{\psi}{\longmapsto} x_i, \\
t_i\stackrel{\psi}{\longmapsto} x^{v_i}t. 
\end{array}
\end{array}
$$
To complete the proof we will show that there is 
an epimorphism $\overline{\psi}$ 
of $K$-algebras that makes this diagram commutative, i.e., 
$\psi=\overline{\psi}\varphi$. To show the existence of
$\overline{\psi}$ we need only show the inclusion
$\ker(\varphi)\subset\ker(\psi)$. As $\ker(\varphi)$, being a toric
ideal, is generated by 
binomials (Sturmfels 1996), it suffices to prove that
any binomial of $\ker(\varphi)$ belongs to $\ker(\psi)$. Let
$$
f=u_0^{a_0}u_1^{a_1}\cdots u_n^{a_n}t_1^{b_1}\cdots t_q^{b_q}-
u_0^{c_0}u_1^{c_1}\cdots u_n^{c_n}t_1^{d_1}\cdots t_q^{d_q}
$$
be a binomial in $\ker(\varphi)$. Then 
$$
t^{a_0}(x_1t)^{a_1}\cdots (x_nt)^{a_n}(x^{v_1}t)^{b_1}\cdots
(x^{v_q}t)^{b_q}=
t^{c_0}(x_1t)^{c_1}\cdots (x_nt)^{c_n}(x^{v_1}t)^{d_1}\cdots
(x^{v_q}t)^{d_q}
$$
Taking degrees in $t$ and $\underline{x}=\{x_1,\ldots,x_n\}$ we obtain
\begin{eqnarray*}
a_0+(a_1+\cdots+a_n)+(b_1+\cdots+b_q)&=&c_0+(c_1+\cdots+c_n)+
(d_1+\cdots+d_q),\\
a_1+\cdots+a_n+2(b_1+\cdots +b_q)&=&c_1+\cdots+c_n+2(d_1+\cdots +d_q).
\end{eqnarray*}
Thus $-a_0+b_1+\cdots+b_q=-c_0+d_1+\cdots+d_q$, and we obtain the
equality
$$
t^{-a_0}x_1^{a_1}\cdots x_n^{a_n}(x^{v_1}t)^{b_1}\cdots
(x^{v_q}t)^{b_q}=
t^{-c_0}x_1^{c_1}\cdots x_n^{c_n}(x^{v_1}t)^{d_1}\cdots
(x^{v_q}t)^{d_q},
$$
i.e., $f\in{\rm ker}(\psi)$, as required. \QED

\medskip

We come to the main result of this section.

\begin{Theorem}\label{round-for-conn-graphs} Let $G$ be a connected
graph 
and let $A$ be its
incidence matrix. Then the system 
$$
x\geq 0;xA\leq\mathbf{1}
$$
has the
integer rounding property if and only if the induced
subgraph of the vertices of any two vertex disjoint odd cycles of $G$ is
connected.
\end{Theorem}

\demo Let $v_1,\ldots,v_q$ be the column vectors of $A$. According 
to (Simis et al. 1998, Theorem~1.1, cf. Villarreal 2005,
Corollary~3.10),  
the subring $K[Gt]:=K[x^{v_1}t,\ldots,x^{v_q}t]$ is normal if and only
if any two vertex disjoint odd cycles of $G$ can be connected by at
least one edge of $G$. Thus we need only show that $K[Gt]$ is normal
if and only the system 
$x\geq 0;xA\leq\mathbf{1}$ has the integer rounding property. 
Let $I=I(G)$ be the edge ideal of $G$.
Since $G$ is connected, by (Simis et al. 1998, Corollary~2.8) the 
subring $K[Gt]$ is normal if and only if the Rees algebra $R[It]$ of
$I$ is normal. By a result of (Herzog et al. 1991), $R[It]$ is normal
if and 
only if $R[It,t^{-1}]$ is normal. By
Proposition~\ref{ext-rees-alg=monsub}, 
$R[It,t^{-1}]$ is normal if and only if the subring 
$$
S=K[t,x_1t,\ldots,x_nt,x^{v_1}t,\ldots,x^{v_q}t]
$$
is normal. Thus we can apply Theorem~\ref{round-up-char} to conclude
that $S$ is normal if and only if the system $x\geq
0;xA\leq\mathbf{1}$ has the integer rounding property.\QED

\begin{Theorem}\label{roundup-iff-rounddown} Let $G$ be a connected
graph and let $A$ be its incidence matrix. Then the system 
$x\geq 0;xA\leq\mathbf{1}$ has the integer rounding property if and
only if any of the following equivalent conditions hold
\begin{description} 
\item{\rm (a)} $x\geq 0;xA\geq\mathbf{1}$ is a system with the
integer rounding property. 
\item{\rm (b)} $R[It]$ is a normal domain, where $I=I(G)$ is the edge
ideal 
of $G$.
\item{\rm (c)} $K[x^{v_1}t,\ldots,x^{v_q}t]$ is normal, 
where $v_1,\ldots,v_q$ are the column vectors of $A$
\item{\rm (d)} $K[t,x_1t,\ldots,x_nt,x^{v_1}t,\ldots,x^{v_q}t]$ is
normal.
\end{description}
\end{Theorem}

\demo According to Theorem~\ref{poset-main-normal}, the system  $x\geq
0;xA\geq\mathbf{1}$ 
has the integer rounding property if and only if the Rees algebra
$R[It]$ 
is normal. Thus the result follows from the proof of
Theorem~\ref{round-for-conn-graphs}. \QED

\section{The canonical module and the $\bf
a$-invariant}\label{section-on-canmod}

In this section we give a description of the canonical module and the
$a$-invariant for subrings arising from systems with the integer
rounding property.

Let $\mathcal{C}$ be a clutter with vertex set $X=\{x_1,\ldots,x_n\}$
and let $v_1,\ldots,v_q$ be the columns of the incidence matrix 
$A$ of $\mathcal{C}$. For use below consider the set $w_1,\ldots,w_r$
of all $\alpha\in\mathbb{N}^n$ such that $\alpha\leq v_i$ for 
some $i$. Let $R=K[x_1,\ldots,x_n]$ be a polynomial ring over a field
$K$ and 
let 
$$
S=K[x^{w_1}t,\ldots,x^{w_r}t]\subset R[t]
$$
be the subring of $R[t]$ generated by
$x^{w_1}t,\ldots,x^{w_r}t$, where $t$ is a new variable. 
As $(w_i,1)$ lies in the hyperplane
$x_{n+1}=1$ 
for all $i$, $S$ is a standard $K$-algebra. Thus a monomial $x^at^b$
in $S$ has degree $b$. In what follows we assume that $S$ has this
grading. Recall that the $a$-{\it invariant\/} of $S$, denoted
$a(S)$, is  
the degree as a rational 
function of the Hilbert series 
of $S$, see for instance (Villarreal 2001, p.~99). If $S$ is 
Cohen-Macaulay and $\omega_S$ is the canonical 
module of $S$, then
\begin{equation}\label{princeton-fall-07}
a(S)=-{\rm min}\{\, i\, \vert\, (\omega_S)_i\neq 0\},
\end{equation}
see (Bruns and Herzog 1997, p.~141) and (Villarreal 2001,
Proposition~4.2.3).  
This formula applies if $S$ is normal because normal 
monomial subrings are Cohen-Macaulay (Hochster 1972). 
If $S$ is normal,
then by a formula of 
Danilov-Stanley, see (Bruns and Herzog 1997, Theorem~6.3.5) and
(Danilov 1978), the canonical module of $S$ is 
the ideal given by  
\begin{equation}\label{ejc8}
\omega_{S}=(\{x^at^b\vert\,
(a,b)\in \mathbb{N}{\mathcal B}\cap ({\mathbb R}_+{\cal B})^{\rm o}\}),
\end{equation}
where  ${\mathcal B}=\{(w_1,1),\ldots,(w_r,1)\}$ and 
$({\mathbb R}_+{\mathcal B})^{\rm o}$ is the interior of 
${\mathbb R}_+{\mathcal B}$ relative to 
${\rm aff}({\mathbb R}_+{\mathcal B})$, 
the affine hull of ${\mathbb R}_+{\mathcal B}$. In our case  
${\rm aff}({\mathbb R}_+{\mathcal B})=\mathbb{R}^{n+1}$. 

The next theorem complements a result of (Brennan et al. 2008). In
loc. cit.  
a somewhat different expression for the canonical module and
$a$-invariant are shown. Our expressions are simpler because 
they only involve the vertices of a certain polytope while in
(Brennan et al. 2008) some 
other parameters are involved.

\begin{Theorem}\label{can-mod-intr} Let $\mathcal{C}$ be a clutter
with incidence matrix 
$A$, let $v_1,\ldots,v_q$ be the columns of $A$, and let 
$w_1,\ldots,w_r$ be the set of all
$\alpha\in\mathbb{N}^n$ such that $0\leq \alpha\leq v_i$ 
for some $i$.
If the system $x\geq 0; xA\leq \mathbf{1}$ has the 
integer rounding property and $\ell_1,\ldots,\ell_m$ are the non-zero
vertices of 
$P=\{x\vert\, x\geq 0;xA\leq\mathbf{1}\}$, then the subring 
$S=K[x^{w_1}t,\ldots,x^{w_r}t]$ is normal, the canonical module of
$S$ is
given by 
\begin{equation}\label{march5-08}
\omega_S=\left(\left\{\left.x^at^b\right\vert\,
(a,b)\left(\begin{array}{rrrlrr} 
-\ell_1&\cdots&-\ell_m&e_1&\cdots& e_n\\
1        &\cdots& 1     &\ 0  &  \cdots&0
\end{array}
\right)>0\right\}\right),
\end{equation}
and the $a$-invariant of $S$ is equal to 
$-(\max_i\{\lfloor|\ell_i|\rfloor\}+1)$. Here
$|\ell_i|=\langle\ell_i,\mathbf{1}\rangle$.
\end{Theorem}

\demo Note that in Eq.~(\ref{march5-08}) we regard the $\ell_i$'s and
$e_i's$ as column vectors. The normality of $S$ follows 
from Theorem~\ref{round-up-char}. 
Let $P=\{x\vert\, x\geq 0;\, xA\leq\mathbf{1}\}$ and let $T(P)$
be its antiblocking polyhedron
$$
T(P):=\{z\vert\, z\geq 0; \langle z,x\rangle\leq 1\mbox{ for all }x\in
P\}. 
$$
By the finite basis theorem (Schrijver 1986) we can write 
\begin{equation}\label{march6-08-1}
P=\{z\vert\, z\geq 0;\langle z,w_i\rangle\leq 1\, \forall i\}={\rm
conv}(\ell_0,\ell_1,\ldots,\ell_m),
\end{equation}
where $\ell_1,\ldots,\ell_m,\ell_0$ are the vertices of $P$  and
$\ell_0=0$. Notice that the vertices of $P$ are in $\mathbb{Q}_+^n$. 
From Eq.~(\ref{march6-08-1}) we
readily get the equality
\begin{equation}\label{march6-08-3}
\{z\vert\, z\geq 0;\langle z,\ell_i\rangle\leq 1\, \forall i\}=T(P).
\end{equation}
Using Eq.~(\ref{march6-08-1}) again and noticing that 
$\langle\ell_i,w_j\rangle\leq 1$ for all $i,j$, we get
$$
\mathbb{R}_+^n\cap({\rm
conv}(\ell_0,\ldots,\ell_m)+\mathbb{R}_+\{-e_1,\ldots,-e_n\})=
\{z\vert\, z\geq 0;\langle z,w_i\rangle\leq 1\, \forall i\}.
$$
Hence using this equality and (Schrijver 1986, Theorem~9.4) we obtain
\begin{equation}\label{march6-08-2}
\mathbb{R}_+^n\cap({\rm
conv}(w_1,\ldots,w_r)+\mathbb{R}_+\{-e_1,\ldots,-e_n\})=
\{z\vert\, z\geq 0;\langle z,\ell_i\rangle\leq 1\, \forall i\}.
\end{equation}
By (Fulkerson 1971, Theorem~8) we have the equality
$$
{\rm conv}(w_1,\ldots,w_r)=\mathbb{R}_+^n\cap({\rm
conv}(w_1,\ldots,w_r)+\mathbb{R}_+\{-e_1,\ldots,-e_n\}).
$$
Therefore using Eqs.~(\ref{march6-08-3})
and  (\ref{march6-08-2}) we conclude the following duality:
\begin{eqnarray}
P=\{x\vert\, x\geq0;\langle x,w_i\rangle\leq{1}\, \forall i\}&=&{\rm 
conv}(\ell_0,\ell_1,\ldots,\ell_m),\nonumber\\
{\rm
conv}(w_1,\ldots,w_r)&=&
\{x\vert\, x\geq0;\langle x,\ell_i\rangle\leq{1}\forall
i\}=T(P).\label{equ-duality}
\end{eqnarray}
We set $\mathcal{B}=\{(w_1,1),\ldots,(w_r,1)\}$. Note that 
$\mathbb{Z}\mathcal{B}=\mathbb{Z}^{n+1}$. 
From Eq.~(\ref{equ-duality}) it is seen that
\begin{equation}\label{march12-08}
\mathbb{R}_+\mathcal{B}=H_{e_1}^+\cap\cdots\cap H_{e_n}^+\cap
H_{(-\ell_1,1)}^+\cap\cdots\cap H_{(-\ell_m,1)}^+.
\end{equation}
Here $H_{a}^+$ denotes the closed halfspace 
$H_a^+=\{x\vert\, \langle
x,a\rangle\geq 0\}$ and $H_a$ stands for the hyperplane through the
origin with normal 
vector $a$. Notice that
$$
H_{e_1}\cap\mathbb{R}_+\mathcal{B},\ldots,H_{e_n}\cap
\mathbb{R}_+\mathcal{B},
H_{(-\ell_1,1)}\cap\mathbb{R}_+\mathcal{B},\ldots,
H_{(-\ell_m,1)}\cap\mathbb{R}_+\mathcal{B}
$$
are proper faces of
$\mathbb{R}_+\mathcal{B}$. Hence from Eq.~(\ref{march12-08}) 
we get that a vector $(a,b)$, with
$a\in\mathbb{R}^n$, $b\in\mathbb{R}$, is in the relative interior of 
$\mathbb{R}_+\mathcal{B}$ if and only if the entries of $a$ are
positive and $\langle(a,b),(-\ell_i,1)\rangle>0$ for all $i$. Thus the
required expression for $\omega_S$, i.e.,
Eq.~(\ref{march5-08}), follows using the normality of $S$
and the Danilov-Stanley formula given in Eq.~(\ref{ejc8}). 

It remains to prove the formula for $a(S)$, 
the $a$-invariant of $S$. Consider the vector 
$(\mathbf{1},b_0)$, where $b_0=\max_i\{\lfloor|\ell_i|\rfloor\}+1$. 
Using Eq.~(\ref{march5-08}), it is not hard to see (by direct
substitution of $(\mathbf{1},b_0)$), that the monomial 
$x^{\mathbf{1}}t^{b_0}$ is in $\omega_S$. Thus from
Eq.~(\ref{princeton-fall-07}) we get $a(S)\geq -b_0$. Conversely if 
the monomial $x^at^b$ is in $\omega_s$, then again from
Eq.~(\ref{march5-08}) we get $\langle(-\ell_i,1),(a,b)
\rangle>0$ for all $i$ and $a_i\geq 1$ for all $i$, where $a=(a_i)$.
Hence 
$$
b>\langle a,\ell_i\rangle\geq \langle
\mathbf{1},\ell_i\rangle=|\ell_i|\geq \lfloor |\ell_i|\rfloor. 
$$
Since $b$ is an integer we obtain $b\geq \lfloor|\ell_i|\rfloor+1$ for
all $i$. Therefore $b\geq b_0$, i.e., $\deg(x^at^b)=b\geq b_0$. As 
$x^at^b$ was an arbitrary monomial in $\omega_s$, by the formula 
for the $a$-invariant of $S$ given in Eq.~(\ref{princeton-fall-07}) 
we obtain that $a(S)\leq -b_0$. Altogether one has $a(S)=-b_0$, as
required.  \QED

\paragraph{Monomial subrings of cliques of perfect graphs}  

Let $S$ be a set of vertices of a graph $G$, 
the {\it induced
subgraph\/} $\langle S\rangle$ is the maximal subgraph of $G$ with
vertex set $S$. A {\it clique\/} of a graph $G$ is a subset of the set
of vertices that induces 
a complete subgraph. Let $G$ be a graph with vertex set
$X=\{x_1,\ldots,x_n\}$. A {\it colouring\/} of the vertices of $G$ is
an assignment 
of colours to the vertices of $G$ in such a way that adjacent vertices
have distinct colours. The {\it chromatic number\/} of $G$ 
is the minimal number of colours in a colouring of $G$.
A graph is {\it perfect\/} if for every induced subgraph $H$, the
chromatic 
number of $H$ equals the size of the largest complete subgraph 
of $H$. Let $S$ be a subset of the vertices of $G$. The set $S$ is called 
{\it independent\/} if no two vertices of $S$ are adjacent. 

For use below we consider the empty set as a clique whose vertex set
is empty. The {\it support\/} of a monomial $x^a$ is given by 
${\rm supp}(x^a)= \{x_i\, |\, a_i>0\}$. Note that ${\rm
supp}(x^a)=\emptyset$ if and only if $a=0$. 

\begin{Theorem}\label{perfect-canon-ainv} Let $G$ be a perfect graph
and let 
$S=K[x^{\omega_1}t,\ldots,x^{\omega_r}t]$ be the subring generated by
all 
square-free monomials $x^at$ such that ${\rm supp}(x^a)$ is a clique
of $G$. Then the canonical module of $S$ is
given by 
$$
\omega_S=\left(\left\{\left.x^at^b\right\vert\,
(a,b)\left(\begin{array}{rrrlrr} 
-a_1&\cdots&-a_m&e_1&\cdots& e_n\\
1        &\cdots& 1       &0  &  \cdots& 0
\end{array}
\right)\geq\mathbf{1}\right\}\right),
$$
where $a_1,\ldots,a_m$ are the characteristic vectors of the
maximal independent sets of $G$, and the $a$-invariant of $S$ is equal to 
$-(\max_i\{|a_i|\}+1)$.
\end{Theorem}

\demo Let $v_1,\ldots,v_q$ be the set of characteristic vectors of the
maximal cliques of $G$. Note that $w_1,\ldots,w_r$ is the set of all
$\alpha\in\mathbb{N}^n$ such that $\alpha\leq v_i$ 
for some $i$. Since $G$ is a perfect graph, by
(Korte and Vygen 2000, Theorem~16.14) we have the equality
$$
P=\{x\vert x\geq 0;xA\leq\mathbf{1}\}={\rm conv}(a_0,a_1,\ldots,a_p),
$$
where $a_0=0$ and $a_1,\ldots,a_p$ are the characteristic vectors of 
the independent sets of $G$. We may assume that $a_1,\ldots,a_m$
correspond to the maximal independent sets of $G$. Furthermore, since
$P$ has only integral vertices, by a result of (Lov\'asz 1972) 
the system $x\geq 0;xA\leq\mathbf{1}$ is totally dual
integral, i.e., the minimum in 
the LP-duality equation
\begin{equation}\label{jun6-2-03}
{\rm max}\{\langle \alpha,x\rangle \vert\, x\geq 0; xA\leq \mathbf{1}\}=
{\rm min}\{\langle y,\mathbf{1}\rangle \vert\, y\geq 0; Ay\geq\alpha\} 
\end{equation}
has an integral optimum solution $y$ for each integral vector $\alpha$ with 
finite minimum. In particular the system $x\geq 0;xA\leq\mathbf{1}$
satisfies the integer rounding property. Therefore the result follows
readily from Theorem~\ref{can-mod-intr}. \QED

\medskip

For use below recall that a graph $G$ is called {\it unmixed\/} if all
maximal independent sets of $G$ have the same cardinality. Unmixed
bipartite graphs have been nicely characterized in (Villarreal 2007).

\begin{Corollary}\label{extended-gorenstein} Let $G$ be a connected
bipartite graph and let 
$I=I(G)$ be its edge ideal. Then the extended Rees algebra 
$R[It,t^{-1}]$ is a Gorenstein standard $K$-algebra if and only if
$G$ is unmixed. 
\end{Corollary}

\demo Let $\omega_S$ be the canonical module of $S=R[It,t^{-1}]$. Recall
that $S$ is Gorenstein if and only if $\omega_S$ is a principal
ideal (Bruns and Herzog 1997). Since any bipartite graph is a perfect
graph, 
the result follows 
using Proposition~\ref{ext-rees-alg=monsub} together with the
description of the canonical module given in
Theorem~\ref{perfect-canon-ainv}. \QED  

\section*{\Large\bf References}

\noindent Brennan JP, Dupont LA and Villarreal RH. 2008. 
Duality, a-invariants and canonical modules         
of rings arising from linear optimization problems. 
Bull. Math. Soc. Sci. Math. Roumanie (N.S.) {51}: 279--305.

\medskip

\noindent {Bruns W and Herzog J. 1997. 
{\em Cohen-Macaulay Rings}. Cambridge University 
Press, Cambridge. Revised Edition.}

\medskip

\noindent {Danilov VI. 1978. The geometry of toric 
varieties. Russian
Math. Surveys {33}: 97-154.}

\medskip

\noindent Dupont LA and Villarreal RH. 2010. Edge ideals of
clique clutters of  
comparability graphs and the normality of monomial ideals.  
Math. Scand. {106}: 88--98.

\medskip

\noindent {Escobar C, Mart\'\i nez-Bernal J and Villarreal RH. 2003. 
Relative volumes and minors
in monomial subrings. Linear Algebra Appl. {374}: 275--290.}

\medskip

\noindent Fulkerson DR. 1971.  Blocking and anti-blocking
pairs of polyhedra.  Math. Programming  {1}: 168--194. 

\medskip

\noindent Gitler I, Reyes E and Villarreal RH. 2009. 
Blowup algebras of square--free monomial ideals and some links to
combinatorial optimization problems. 
Rocky Mountain J. Math. {39}: 71--102. 

\medskip

\noindent {Herzog J, Simis A and  Vasconcelos WV. 1991. 
Arithmetic of
normal Rees algebras. J. Algebra {143}: 269--294.} 

\medskip

\noindent {Hochster M. 1972. Rings of invariants of tori,
{C}ohen-{M}acaulay rings generated by
 monomials, and polytopes. Ann. of Math. 
{96}: 318--337.}

\medskip

\noindent Korte B and  Vygen J. 2000. 
{\it Combinatorial Optimization Theory 
and Algorithms}. Springer-Verlag.

\medskip

\noindent Lov\'asz L. 1972. Normal hypergraphs and the perfect
graph conjecture. Discrete Math. {2}: 253--267. 

\medskip

\noindent {Simis A, Vasconcelos WV and Villarreal RH. 1994. On the
ideal theory of graphs. J. Algebra {167}: 389--416.}

\medskip

\noindent {Simis A, Vasconcelos WV and Villarreal RH. 1998. The
integral closure of subrings associated to graphs. J. Algebra {199}:
281--289.}   

\medskip

\noindent {Schrijver A. 1986. {\it Theory of Linear and 
Integer Programming}. John Wiley \& Sons, New York.}

\medskip

\noindent {Sturmfels B. 1996. {\em Gr\"obner Bases and Convex
Polytopes\/}. University Lecture Series {8}. American 
Mathematical Society, Rhode Island.}  

\medskip

\noindent Valencia C and Villarreal RH. 2003. Canonical
modules of certain edge subrings. 
European J. Combin. {24}: 471--487.

\medskip

\noindent Vasconcelos WV. 2005. {\it Integral Closure}. Springer
Monographs in Mathematics, Springer, New York.

\medskip

\noindent {Villarreal RH. 2001. {\it Monomial Algebras}. Monographs
and  
Textbooks in Pure and Applied Mathematics {238}. Marcel 
Dekker, New York.} 

\medskip

\noindent Villarreal RH. 2005. Normality of semigroups with some 
links to graph theory. Discrete Math. {302}: 267-284.

\medskip

\noindent Villarreal RH. 2007. Unmixed bipartite graphs. Rev.
Colombiana Mat. {41}: 393--395.

\end{document}